\documentclass[secthm, 10pt]{elsart}

\usepackage{amsmath}
\usepackage{amsopn}
\usepackage{amssymb}
\usepackage[all]{xy}

\newcommand{\br}[1]{\overline{#1}}

\newcommand{\td}[1]{\widetilde{#1}}

\newcommand{\ZZ}{\mathbb{Z}}
\newcommand{\RR}{\mathbb{R}}

\newcommand{\QQ}{\mathbb{Q}}
\newcommand{\WW}{\mathbb{W}}
\newcommand{\FF}{\mathbb{F}}

\newcommand{\MS}{\mathbb{S}}
\newcommand{\MSL}{{\mathbb{S}l}}
\newcommand{\MSy}{{\mathbb{S}}^0}
\newcommand{\MSyL}{{\mathbb{S}l}^0}

\numberwithin{equation}{section}
\numberwithin{figure}{section}
\DeclareMathOperator{\Ext}{Ext}
\DeclareMathOperator{\Hom}{Hom}
\DeclareMathOperator{\End}{End}
\DeclareMathOperator{\aut}{Aut}

\DeclareMathOperator{\Spec}{Spec}

\newcommand{\Tr}{\mathop{Tr}}

\begin{document}

\begin{frontmatter}

\title{Isogenies of elliptic curves and the Morava stabilizer group}
\author{Mark Behrens\thanksref{NSF}}
\address{
Department of Mathematics,
MIT,
Cambridge, MA 02139, USA}
\author{Tyler Lawson\thanksref{NSF}}
\address{
Department of Mathematics,
MIT,
Cambridge, MA 02139, USA}
\thanks[NSF]{The authors were supported by the NSF.}

\begin{abstract}
Let $\MS_2$ be the $p$-primary second Morava stabilizer group, $C$ a
supersingular elliptic curve over $\br{\FF}_{p}$, $\mathcal{O}$ the ring of
endomorphisms of $C$,
and $\ell$ a topological generator of
$\ZZ_p^\times$ (respectively $\ZZ_2^\times/\{\pm 1\}$ if $p=2$).  
We show that for $p > 2$ the group $\Gamma \subseteq 
\mathcal{O}[1/\ell]^\times$ 
of quasi-endomorphisms of degree a power of $\ell$ is dense in $\MS_2$.
For $p=2$, we show that $\Gamma$ is dense in an index $2$ subgroup of
$\MS_2$.

\end{abstract}


\end{frontmatter}

\section*{Introduction}

Fix a prime $p$.  Let $F_n$ be the Honda height $n$ formal group over
$\FF_{p^n}$.  The endomorphism ring $\mathcal{O}_p = \End(F_n)$ 
is the unique maximal order of the
$\QQ_p$-division algebra $D_p$ of Hasse invariant $1/n$
\cite{Ravenel}, and the Morava stabilizer group $\MS_n$ is the
automorphism group $\aut(F_n) = \mathcal{O}_p^\times$.  This group is a
$p$-adic analytic group of dimension $n^2$, and is of interest to
topologists because it is intimately related to the $n$th layer of the
chromatic filtration on the stable homotopy groups of spheres.
We wish to understand the group $\MS_n$ for $n=2$
from the point
of view of elliptic curves.

Throughout this paper we let $C$ be a fixed supersingular elliptic
curve defined over $\FF_{p^2}$.  Let $\mathcal{O} = \End(C)$ be the
ring of endomorphisms of the curve $C$ defined over $\br{\FF}_p$, and let 
$D = \mathcal{O}\otimes \QQ$ be the ring of
quasi-endomorphisms.  Because $C$ is supersingular, it is known
that $D$ is the quaternion algebra over $\QQ$ ramified at $p$
and $\infty$ \cite{Kohel}, and that $\mathcal{O}$ is a maximal order of
$D$ \cite[3.1]{Silverman}, \cite{Deuring}.  The reduced norm
$$ N: D \rightarrow \QQ_p^\times $$
gives the degree of the quasi-endomorphism.
Let $\widehat{C}$ be the formal
completion of $C$ at the identity.  Because $C$ is supersingular, the formal
group $\widehat{C}$ is isomorphic to the Honda formal group $F_2$ over
$\br{\FF}_p$.  In
fact, Tate proved that the natural map
$$ \rho : \End(C)\otimes \ZZ_p  \rightarrow \End(\widehat{C}) $$
is an isomorphism 
\cite{Waterhouse-Milne}.
The isomorphism $\rho$ extends to an isomorphism
$$ \rho' : D \otimes_\QQ \QQ_p \xrightarrow{\cong} D_p, $$
making explicit the fact that $D$ is ramified at $p$.  The map $\rho'$ is
compatible with the reduced norm map on the division algebras $D$ and
$D_p$.

Fix $\ell \ge 2$ to be coprime to $p$.  As the notation
suggests, we intend for $\ell$ to be another prime, but this is
unnecessary for the results of this paper.
Define a monoid 
$$ \Gamma = \{x \in \mathcal{O}[1/\ell] \: : \:
N(x) \in \ell^\ZZ \} \subseteq \mathcal{O}[1/\ell]^\times .$$
Then $\Gamma$ is actually a group: given an endomorphism
$\phi$ of degree $\ell^k$, the quasi-endomorphism $\ell^{-k}\widehat{\phi}$
is its inverse, where $\widehat{\phi}$ is the dual isogeny.
Note that if $\ell$ is prime, then $\Gamma = \mathcal{O}[1/\ell]^\times$.
The group $\Gamma$ may be regarded as being contained in the group $\MS_2$
using the map $\rho$.
The purpose of this note is to prove the following
theorem.

\begin{thm}\label{thm:Gamma}
Suppose that $\ell$ is a topological generator of the group $\ZZ_p^\times$
(respectively the group $\ZZ_2^\times/\{\pm 1 \}$ for $p = 2$).  For $p
> 2$, the group $\Gamma$ is dense
in $\MS_2$.  For $p = 2$, the group $\Gamma$ is dense in the index $2$
subgroup $\td{\MS}_2$, which is the kernel of the composite
$$ \MS_2 \xrightarrow{N} \ZZ_2^\times \rightarrow
(\ZZ/8)^\times/\{1,\ell\}. $$
\end{thm}

Let $\MSL_2$ be the kernel of the reduced norm, so that there is an
exact sequence
$$ 1 \rightarrow \MSL_2 \rightarrow \MS_2 \xrightarrow{N} \ZZ_p^\times
\rightarrow 1.
$$
Similarly, let $\Gamma^1$ be the corresponding subgroup of $\Gamma$,
so that there is an exact sequence
$$ 1 \rightarrow \Gamma^1 \rightarrow \Gamma \xrightarrow{N} \ell^\ZZ
\rightarrow 1.
$$
(We will see in the proof of Theorem~\ref{thm:Gamma} that $N$ is indeed
surjective.)
We denote by $\MSy_2$ the $p$-Sylow subgroup of $\MS_2$, so that there is a
short exact sequence
$$ 1 \rightarrow \MSy_2 \rightarrow \MS_2 \rightarrow \FF_{p^2-1}^\times
\rightarrow 1. $$
Similarly we let $\MSyL_2$ be the subgroup of $\MSy_2$ of elements of norm
$1$.  Define $\Lambda$ to be the subgroup $\Gamma^1 \cap \MSyL_2$.
Theorem~\ref{thm:Gamma} will follow from the following norm $1$ versions.
Note that in the following theorem and corollary, 
$\ell$ is only assumed to be relatively prime
to $p$.

\begin{thm}\label{thm:Lambda}
The group $\Lambda$ is dense in $\MSyL_2$.
\end{thm}

\begin{cor}\label{cor:Gamma1}
The group $\Gamma^1$ is dense in $\MSL_2$.
\end{cor}

We pause to explain the reason why these theorems are interesting from
the point of view of homotopy theory.  
The $p$-component of the stable homotopy groups of spheres admits an
especially rich filtration known as the chromatic filtration
\cite{Ravenelorange}.
Work of Morava, Hopkins, Miller, Goerss, and Devinatz \cite{Morava},
\cite{RezkHMT}, \cite{GoerssHopkins}, \cite{DevinatzHopkins} shows that 
the group $\MS_n$ acts on the Morava $E$-theory spectrum $E_n$, and
the $n$th layer of the chromatic filtration is described by the homotopy
fixed points $E_n^{h\MS_n}$ of this action.  The first chromatic layer is
completely understood.  The second chromatic layer is currently the subject of
intense study.  Goerss, Henn, Mahowald, and Rezk \cite{GHMR} produced a
decomposition of $E_2^{h\MS_2}$ at the prime $3$ 
in terms of finite homotopy fixed point
spectra. The first author gave an interpretation of their work 
in terms of the moduli
space of elliptic curves in \cite{BehrensModular}.  In that paper, a
spectrum $Q(\ell)$ was introduced which was a shown to be a good
approximation to $E_2^{h\MS_2}$ for $p=3$ and $\ell = 2$.  
In future work \cite{BehrensTree}, we
will show that the spectrum $Q(\ell)$
is the homotopy fixed point spectrum $E^{h\Gamma}$.  In particular,
Theorem~\ref{thm:Gamma} shows that, in some sense, the spectrum 
$Q(\ell)$ is a good
approximation for $E_2^{h\MS_2}$ for all $p$ and suitable $\ell$.

Gorbounov, Mahowald, and Symonds \cite{GMS} studied dense subgroups of
$\MSyL_n$ ($\MS l$ in their notation), and we prove 
Theorem~\ref{thm:Lambda} using their methods. 
In particular, it is shown in \cite{GMS} that if $p = 3$, 
then there is a dense subgroup $\ZZ/3 \ast \ZZ/3$
contained in $\MSyL_2$.  
In \cite{BehrensTree}, it will be shown that for
$\ell=2$, the group $\Lambda$ is $\ZZ/3 \ast \ZZ/3$.
More generally, the groups $\Lambda$ and $\Gamma^1$ admit
explicit presentations as finite amalgamations for any $p$ and $\ell$.

The authors were alerted by the referee to the related work of Baker
\cite{Baker}.  Baker studies, for primes $p \ge 5$, the category whose
objects are supersingular elliptic curves over $\br{\FF}_p$, and whose
morphisms are the morphisms of the associated formal groups.  He then
proves an analog of Morava's change of rings theorem: roughly
speaking, he shows that the continuous cohomology of this category of
supersingular curves computes the $E_2$-term of the $K(2)$-local
Adams-Novikov spectral sequence converging to $\pi_*(S_{K(2)})$.  The
chief difference between this paper and the work of Baker is that we
insist on working only with the actual rings of isogenies, and not
their $p$-completions.

In Section~\ref{sec:super}, we recall the relationship between maximal
orders of $D$ and the endomorphism rings of supersingular curves at $p$.  In
Section~\ref{sec:density}, we recall the homological criterion that is
employed in \cite{GMS} to detect dense subgroups of $\MSyL_2$ for $p >2$.
We then extend these methods to give an explicit criterion for density at
the prime $2$.
We use
these criteria in Section~\ref{sec:proofs} to prove
Theorem~\ref{thm:Gamma}, Theorem~\ref{thm:Lambda}, and 
Corollary~\ref{cor:Gamma1}.

The first author would like to thank Hans-Werner Henn for pointing
out to him that Theorem~\ref{thm:Gamma} was true in the case of $p =
3$ and $\ell=2$.  Thanks also go to Johan de Jong and Catherine O'Neil
for helpful discussions related to this paper.

\section{Supersingular curves and endomorphism rings}\label{sec:super}

In this section we recall the correspondence between endomorphism rings of
supersingular curves and maximal orders of $D$. 
Any two supersingular elliptic curves $C_1$ and $C_2$ over $\br{\FF}_{p}$ are 
isogenous.  In fact, Kohel proves the following proposition.

\begin{prop}[Kohel {\cite[Cor. 77]{Kohelthesis}}]\label{prop:Kohel}
Let $C_1$ and $C_2$ be supersingular elliptic curves over $\br{\FF}_{p}$.
Then for all $k \gg 0$, there exists an isogeny $\phi: C_1 \rightarrow C_2$
of degree $\ell^k$.
\end{prop}

Let $X^{ss}$ be the collection of isomorphism classes of supersingular
curves $C'$ over $\br{\FF}_p$.  Given an isogeny $\phi: C \rightarrow C'$
of degree $N$, we define a map
$$ \iota_\phi: \End(C') \rightarrow \End(C)\otimes \ZZ[1/N] \subset D$$
by 
$$\iota_\phi(\alpha) = \frac{1}{N} \cdot \widehat{\phi} \circ \alpha
\circ \phi.$$
The map $\iota_\phi$ is a ring homomorphism and its image is a maximal order in
$D$.  If $\phi':C \rightarrow C'$ is another choice of isogeny then it is
easily seen that the maximal order $\iota_{\phi'}(\End(C'))$ is conjugate
to $\iota_{\phi}(\End(C'))$.  Let $\mathcal{M}_{D}$ be the collection of
conjugacy classes of maximal orders of $D$.  Consider the map
$$ \xi: X^{ss} \rightarrow \mathcal{M}_D $$
given by $\xi([C']) = [\iota_{\phi}(\End(C'))]$.  

\begin{thm}[Deuring \cite{Deuring}, Kohel \cite{Kohel}]\label{thm:Deuring}
The map $\xi$ is a surjection and the preimage of a conjugacy class 
$[{\mathcal O}']$ of maximal orders consists of either a single class
represented by a curve with $j$-invariant in $\FF_p$, or two classes
represented by curves with distinct Galois-conjugate $j$-invariants in
$\FF_{p^2}$.
\end{thm}

In the former case, the elliptic curve can be defined over $\FF_p$,
and in the latter case it can only be defined over $\FF_{p^2}$.

We choose a preferred set of representatives of conjugacy classes in
$\mathcal{M}_D$ for the remainder of this note.
Fix a choice of representative $C'$ of each
isomorphism class $[C'] \in X^{ss}$.
Using Proposition~\ref{prop:Kohel}, choose for each $C'$ an isogeny
$$ \phi_{C'}: C \rightarrow C' $$
of degree $\ell^{e(C')}$.  
Define $\mathcal{O}_{C'}$ to be the maximal order
$\iota_{\phi_{C'}}(\End(C'))$.
By letting $C$ be the representative of its isomorphism class, and fixing 
$\phi_C = \mathrm{Id}_C$, we can arrange that $\mathcal{O}_{C} = \End(C)$.
The following is immediate from
Theorem~\ref{thm:Deuring}.

\begin{cor}\label{cor:orders}
Every maximal order $\mathcal{O}'$ of $D$ is conjugate to one of the form
$\mathcal{O}_{C'}$ for some $C' \in X^{ss}$.
\end{cor}

Let $\mathcal{O}'$ be a maximal order of $D$.  Then by
Corollary~\ref{cor:orders}, $c_y(\mathcal{O}')
= y^{-1} \mathcal{O}' y$ is equal to $\mathcal{O}_{C'}$ for some $y
\in D$ and some supersingular elliptic curve $C'$.
The map $\iota_{\widehat \phi} \circ c_y : \mathcal{O}' \rightarrow
\End(C) \otimes \ZZ[1/\ell]$ extends to a norm-preserving ring isomorphism
$$ \iota_{\widehat \phi} \circ c_y : \mathcal{O}' \otimes
\ZZ[1/\ell] \to \mathcal{O} \otimes \ZZ[1/\ell]$$
with inverse $c_{y^{-1}} \circ \iota_\phi$.  In
particular, we have the following.

\begin{cor}\label{cor:minpoly}
Suppose that $x'$ is contained in a maximal $\ZZ[1/\ell]$-order
$\mathcal{O}'[1/\ell]$ of $D$.  Then there exists an element $x \in
\mathcal{O}[1/\ell]$ with the same minimal polynomial as $x'$.
\end{cor}



\section{A cohomological criterion for density}\label{sec:density}

In this section we recall some material from \cite{GMS}, but we give this
material a slightly different treatment.  Our reason is that the 
authors of \cite{GMS} use results of Riehm \cite{Riehm} on the structure of
the commutator subgroups of $\MSL_n$.  Riehm's analysis, however,
excludes the case of $n=2$ and $p=2$, and it turns out that this case has
different behavior.

The maximal order
$\mathcal{O}_p$ of $D_p$ admits a presentation \cite[Appendix~2]{Ravenel}
\begin{equation}\label{eq:Op}
\mathcal{O}_p = \WW\langle S\rangle/ (S^2 = p, Sa = \br{a}S).
\end{equation}
Here $\WW = \WW(\FF_{p^2})$ is the Witt ring with residue field
$\FF_{p^2}$, and $\br{a}$ denotes the Galois conjugate
(lift of the Frobenius on $\FF_{p^2}$) of an element $a \in \WW$.
Every element of
$\mathcal{O}_p$ can then be written uniquely in the form
$$ a + b S $$
for $a,b \in \WW$.  The group $\MS_2 = \mathcal{O}_p^\times$ 
consists of all such elements where $a \not \equiv 0 \pmod p$.

Let $\MSy_2$ be the $p$-Sylow subgroup of $\MS_2$.  The group $\MSy_2$
consists of all elements $a + bS$ where $a \equiv 1 \pmod p$.  
The subgroup $\MSyL_2$ of elements of 
$\MSy_2$ of norm $1$
is the $p$-Sylow subgroup of $\MSL_2$.

Suppose that $G$ is a pro-$p$-group.
Let $G^*$ be the Frattini subgroup of $G$, which is the minimal closed
normal subgroup that contains $G^p$ and $[G,G]$.
Then we have the following theorem.

\begin{thm}[Koch \cite{Koch}, Serre \cite{Serre}]
Suppose that $H$ is a subgroup of a pro-$p$-group $G$.  
Then $H$ is dense in $G$ if
and only if the composite
$$ H \hookrightarrow G \rightarrow G/G^* = H_1^c(G; \FF_p) $$
is surjective.
\end{thm}

\begin{cor}\label{cor:H1}
Let $\{\alpha_i\}$ form an $\FF_p$-basis of the continuous group
homomorphisms $\Hom^c(G,\FF_p) = H^1_c(G;\FF_p)$.  Then $H$ is a dense
subgroup of $G$ if and only if the composite
$$ H \rightarrow G \xrightarrow{\bigoplus_i \alpha_i} 
\bigoplus_i \FF_p $$
is surjective.
\end{cor}

Every
element $x \in \MSy_2$ may be written uniquely in the form 
\begin{equation}\label{eq:ti}
x = (1 + pt_2 + p^2t_4 + \cdots) + (t_1 + pt_3 + p^2t_5 + \cdots)S,
\end{equation}
where $t_i = t_i(x)$ are Teichm\"uller lifts of elements of 
$\FF_{p^2}$ in $\WW$.  This is equivalent to saying that the elements $t_i$
satisfy $t_i^{p^2} = t_i$.  The coefficients $t_i$ give rise to continuous 
functions
$$ t_i: \MSy_2 \rightarrow \FF_{p^2}. $$
Ravenel \cite{RavenelStructure} uses this presentation to express $\MSy_2$ as
the $\FF_{p^2}$-points of a pro-affine group scheme $\Spec S(2)$, where $S(2)$
is the Morava stabilizer algebra
$$ S(2) = \FF_p[t_1, t_2, t_3, \ldots ]/(t_i^{p^2}=t_i). $$
The algebra $S(2)$ is a Hopf algebra.
Equation (\ref{eq:ti}) gives an isomorphism of groups
$$ \MSy_2 \cong \Spec S(2)(\FF_{p^2}). $$
Ravenel shows that this isomorphism gives an isomorphism in cohomology:
\begin{align*}
H^*_c(\MSy_2; \FF_{p^2}) & = H^*(S(2)) \otimes_{\FF_p} \FF_{p^2} \\
& = \Ext_{S(2)}(\FF_p,\FF_{p}) \otimes_{\FF_p} \FF_{p^2}.
\end{align*}
The Ext group is taken in the category of $S(2)$-comodules.

For an arbitrary element $x \in \MSy_2$ expressed as in equation~(\ref{eq:ti}),
express the norm $N(x) \in \ZZ_p^\times$ by 
$$ N(x) = 1 + ps_1 + p^2s_2 + \cdots,$$
where the elements $s_i = s_i(t_1,t_2, \ldots)$ are polynomial functions of
the $t_i$, and $s_i^p = s_i$ are Teichm\"uller lifts of elements of
$\FF_{p}$ (compare with the discussion preceding
Theorem~6.3.12 of \cite{Ravenel}).
Then we may define a quotient Hopf algebra
$$ Sl(2) = S(2)/(s_i(t_1, t_2, \ldots)), $$
whose $\FF_{p^2}$ points give the subgroups $\MSyL_2$ and for which
$$ H^*_c(\MSyL_2; \FF_{p^2}) = H^*(Sl(2)) \otimes_{\FF_p} \FF_{p^2}. $$

The following computation is obtained from combining Theorems~6.2.7 and
6.3.12 of \cite{Ravenel}.

\begin{lem}
For $p > 2$,  we have
$$ H^1(Sl(2)) = \FF_p\{ h_{1,0}, h_{1,1} \}, $$
where $h_{1,i}$ is represented by the element $[t_1^{p^i}]$ in the cobar
complex for $Sl(2)$.
\end{lem}

\begin{cor}
For $p > 2$, the group 
$$ H^1_c(\MSyL_2;\FF_{p^2}) \cong \Hom^c(\MSyL_2,\FF_{p^2}) $$
has an $\FF_{p^2}$-basis consisting of the continuous homomorphisms
$$ t_1, t_1^p : \MSy_2 \rightarrow \FF_{p^2}. $$
\end{cor}

\begin{cor}[Gorbounov-Mahowald-Symonds \cite{GMS}]\label{cor:density}
For $p > 2$, a subgroup $H$ of $\MSyL_2$ is dense if and only if the 
composite
$$ H \hookrightarrow \MSy_2 \xrightarrow{t_1} \FF_{p^2} $$
is surjective.
\end{cor}

\begin{pf*}{Proof.}
Let $\omega \in \FF_{p^2}$ be a primitive $p^2-1$ root of unity.  We may
compute the cohomology with $\FF_{p}$ coefficients by taking $Gal =
Gal(\FF_{p^2}/\FF_p)$ fixed points, and obtain
$$ \Hom^c(\MSy_2,\FF_p) = \Hom^c(\MSy_2,\FF_{p^2})^{Gal}. $$
(The Galois group only acts on the coefficient group and not on $\MSy_2$.)
Here the Frobenius $\sigma \in Gal$ acts by $\sigma(t_1^{p^i}) =
t_1^{p^{i+1}}$
for $i \in \ZZ/2$.  An $\FF_p$-basis for this fixed-point module is given
by the pair of homomorphisms
$$ t_1+t_1^p, \omega t_1 + \omega^p t_1^p: \MSy_2 \rightarrow \FF_{p^2}. $$
The result now follows from Corollary~\ref{cor:H1}.
\qed \end{pf*}

We now address the case where $p=2$.

\begin{lem}
Let $p=2$. Then we have
$$ H^1(Sl(2)) = \FF_2\{ h_{1,0}, h_{1,1}, h_{3,0}, h_{3,1} \},$$
where the generators are represented in the cobar
complex for $Sl(2)$ by
\begin{align*}
h_{1,i} & = [t_1^{2^i}], \\
h_{3,i} & = [(t_3 + t_1t_2)^{2^i}].
\end{align*}
\end{lem}

\begin{pf*}{Proof.}
We follow the same approach of \cite[6.3]{Ravenel}
using the May spectral sequence.  (It is important to refer to the
second edition of \cite{Ravenel}; the previous
version, as well as \cite{RavenelCohomology}, had an error in the restriction
formula in the restricted Lie algebras $\td{L}(n)$.)
The May spectral sequence for $Sl(2)$ takes the form
$$ E_2^{s,*} = H^s(E^0Sl(2)) \Rightarrow H^s(Sl(2)).  $$
The $E_1$-term may be regarded as the Koszul complex for $(E^0Sl(2))^*$
$$ E_1^{*,*} = \FF_2[h_{i,j} \: : \: i \ge 1, j \in \ZZ/2]/(h_{2k,j} +
h_{2k,j+1}), $$
with differential
$$ d_1(h_{i,j}) = 
\begin{cases}
\sum_{i_1+i_2 = i} h_{i_1,j}h_{i_2,j+i_1} & i \le 4, \\
h_{i-2,j+1}^2 & i > 4. 
\end{cases} $$
We see that the only elements of $E_1^{1,*}$ that persist to $E_2^{1,*}$ are
$h_{1,0}$, $h_{1,1}$, $h_{3,0}$, and $h_{3,1}$.  

We will show that these elements are permanent cycles in the May spectral
sequence by explicitly producing cocycles in the cobar complex that they
detect.
By taking the images of the formulas for the coproduct on $BP_*BP$ in 
\cite{Giambalvo}, we arrive at the following formulas for the coproduct in 
$Sl(2)$.
\begin{align*}
\Delta(t_1) & = t_1 \otimes 1 + 1 \otimes t_1, \\
\Delta(t_2) & = t_2 \otimes 1 + t_1 \otimes t_1^2 + 1 \otimes t_2, \\
\Delta(t_3) & = t_3 \otimes 1 + t_1 \otimes t_2^2 + t_2 \otimes t_1 + 
t_1^2 \otimes t_1^2 + 1 \otimes t_3.
\end{align*}
Using the relation
$$ s_1 = t_2 + t_2^2 + t_1^3 = 0, $$
these formulas may be used to verify that the cobar expressions in the 
statement of the lemma are permanent cycles.
\qed \end{pf*}

\begin{cor}
For $p=2$, the group 
$$ H^1_c(\MSyL_2;\FF_4) \cong \Hom^c(\MSyL_2,\FF_4) $$
has an $\FF_4$-basis given by the continuous homomorphisms
$$ t_1, t_1^2, t_3+t_1t_2,(t_3+t_1t_2)^2 : \MSyL_2 \rightarrow \FF_4. $$
\end{cor}

\begin{cor}\label{cor:density_2}
For $p=2$, a subgroup $H$ of $\MSyL_2$ is dense if and only if the composite
$$ H \hookrightarrow \MSy_2 \xrightarrow{t_1 \oplus (t_3+t_1t_2)} 
\FF_4\oplus \FF_4 $$
is surjective.
\end{cor}

\begin{pf*}{Proof.}
Let $\omega \in \FF_4$ be a primitive $3$rd root of unity.  Just as in
Corollary~\ref{cor:density}, we
compute the cohomology with $\FF_{2}$ coefficients by taking $Gal =
Gal(\FF_{4}/\FF_2)$ fixed points, and obtain
$$ \Hom^c(\MSy_2,\FF_2) = \Hom^c(\MSy_2,\FF_{4})^{Gal}. $$
(As before, the Galois group only acts on the coefficient group.)
An $\FF_2$-basis for this fixed-point module is given
by the homomorphisms $t_1+t_1^2$, $\omega t_1 + \omega^2 t_1^2$,
$t_3 + t_1t_2 + (t_3 + t_1t_2)^2$, and $\omega(t_3 + t_1t_2) +
\omega^2(t_3+t_1t_2)^2$.
\qed \end{pf*}

\section{Proof of Theorem~\ref{thm:Gamma}, Theorem~\ref{thm:Lambda}, and
Corollary~\ref{cor:Gamma1}}\label{sec:proofs}

We will make use of the following proposition, which is a special case of
Proposition~9.19 of \cite{Swan}.

\begin{prop}\label{prop:Swan}
Suppose that $f(x) = x^2 + a_1 x + a_2$ is a monic polynomial over $\QQ$
that is irreducible over $\QQ_p$ and $\RR$.  Then there exists an $\alpha$
in $D$ with $f(\alpha) = 0$.  If the elements $a_i$ are integral over
$R \subset \QQ$, then $\alpha$ lies in a maximal $R$-order of $D$. 
\end{prop}

\begin{pf*}{Proof of Theorem~\ref{thm:Lambda} for $p > 2$.}
We will use Proposition~\ref{prop:Swan} to produce elements $x_1,x_2$ in
$\Lambda$ so that $t_1(x_1)$, $t_1(x_2)$ form an $\FF_p$ basis of
$\FF_{p^2}$. Corollary~\ref{cor:density} then yields the result.

Choose integers $r_1$ and $r_2$ such that $r_i \not \equiv 0 \pmod p$, and
so that $r_1$ is a square and $r_2$ is not a square in $\FF_p$.  Let
$$ \alpha_i = \frac{-pr_i - 2}{\ell^{m_ip(p-1)}}, $$
where the integers $m_i$ are chosen sufficiently large so that
\begin{equation}\label{eq:infty}
\alpha_i^2 < 4.
\end{equation}
We claim that the polynomials
$$ f_i(x) = x^2 + \alpha_i x + 1 $$
are irreducible over $\RR$ and $\QQ_p$.  It suffices to check that the
discriminants $\Delta_i = \alpha_i^2 - 4$ are not squares in each of these
fields.  Condition~(\ref{eq:infty}) guarantees that $\Delta_i$ is not a
square in $\RR$.  Over $\QQ_p$ we note that $\Delta_i$ lies in $\ZZ_p$, so
it suffices to check that $\Delta_i$ is not a square in $\ZZ/p^2$.
Because $\ell^{p(p-1)}$ is congruent to $1$ in $\ZZ/p^2$, we have
$$ \Delta_i \equiv 4pr_i \pmod {p^2}. $$
As $r_i$ is not congruent to $0 \pmod p$, $\Delta_i$ is not a square in
$\ZZ/p^2$.

Applying Proposition~\ref{prop:Swan}, we see that there exist $\td{x}_i$
in $D$ so that $f_i(\td{x}_i) = 0$.  
The elements $\td{x}_i$ satisfy monic quadratics over
$\ZZ[1/\ell]$, and so these elements are contained in maximal
$\ZZ[1/\ell]$ orders $\mathcal{O}_i[1/\ell]$ of $D$.  
Applying Corollary~\ref{cor:minpoly}, there exist elements $x_i \in
\mathcal{O}[1/\ell]$ such that 
\begin{equation}\label{eq:chareq}
f_i(x_i) = x_i^2 + \alpha_i x_i + 1 = 0.
\end{equation}
The $x_i$ satisfy $N(x_i) = 1$.  Therefore, we conclude that the
elements $x_i$ are contained in the group 
$\Gamma^1$.

The images of the $\alpha_i$ in $\QQ_p$ lie in $\ZZ_p$, so the 
images of the $x_i$ in $D_p$ lie in $\mathcal{O}_p$. 
Write $x_i$ in the form
$$ x_i = a_i + b_iS $$ 
for $a_i,b_i \in \WW$. 
Reducing equation~(\ref{eq:chareq}) modulo the ideal $(S)$, we see
that
$$ x_i^2 - 2x_i +1 \equiv a_i^2 - 2 a_i + 1 \equiv 0 \pmod S. $$
We conclude that $a_i \equiv 1 \pmod p$.
This implies that the elements $x_i$ actually lie in $\Lambda$, and their
images in $\mathcal{O}_p$ are of the form
$$ x_i = (1 + pa_i' + b_i S) $$
for $a_i' \in \WW$.  

Equation~(\ref{eq:chareq}) implies that the reduced trace of $x_i$ is
given by
\begin{equation}\label{eq:trace}
\Tr(x_i) = 2 + p\Tr(a_i') = -\alpha_i,
\end{equation}
whereas the reduced norm is given by
\begin{equation}\label{eq:norm}
N(x_i) = 1 + p\Tr(a_i') + p^2N(a_i') - pN(b_i) = 1.
\end{equation}
Substituting the expression for $\Tr(a_i')$ given by
equation~(\ref{eq:trace}) gives
\begin{equation}\label{eq:master}
N(b_i) = pN(a_i') - \frac{\alpha_i+2}{p}.
\end{equation}
Note that $\alpha_i \equiv -2 \pmod p$.
Reducing equation~(\ref{eq:master}) modulo $p$ yields
$$ N(t_1(x_i)) \equiv N(b_i) \equiv r_i \in \FF_p. $$
If $t_1(x_1)$ and $t_1(x_2)$ were $\FF_p$ linearly dependent in
$\FF_{p^2}$,
their norms would lie in the same quadratic residue class.  The
$r_i$ were chosen so that this does not happen, so we conclude that
$\{t_1(x_1), t_1(x_2)\}$ forms a basis of $\FF_{p^2}$.  Therefore, by
Corollary~\ref{cor:density}, the subgroup $\Lambda$ is dense in
$\MSyL_2$.
\qed \end{pf*}

\begin{pf*}{Proof of Theorem~\ref{thm:Lambda} for $p = 2$.}
The proof is similar to the proof for $p > 2$, but more involved.
There is precisely one isomorphism class of supersingular elliptic
curve $C$ at $p = 2$.  It follows that $D$ has one conjugacy class
of maximal order.  By checking the invariants of the division
algebra $D$, it can be shown \cite{Pizer} that $D$ is of the form of
the rational quaternions
$$ \QQ\langle i , j \rangle /(i^2 = j^2 = -1, ij = -ji).$$
We may therefore assume that $\End(C) \subset D$ is the maximal order
$\mathcal{O}$
generated by
$$ \{\omega, i, j, k\}, $$
where $k = ij$ and $\omega = \frac{1 + i + j + k}{2}$.  Note that $\omega^3
= 1$.  The automorphism group $\aut C = \End(C)^\times$ is the binary
tetrahedral group $\tilde A_4$ of order $24$ given by the semidirect
product $Q_8 \rtimes C_3$.  The cyclic 
group $C_3$ is generated by $\omega$ and the quaternion group $Q_8$ is
generated by $i$ and $j$.  We have $\omega i \omega^2 = j$ and $\omega j
\omega^2 = k$.

Let $T$ be the element $i - j \in \mathcal{O}$.  Then we have $T^2 = -2$
and $T \omega = \omega^2 T$.  The Witt ring $\WW = \WW(\FF_4)$ will be
identified with the subring
$$ \ZZ_2[\omega] \subset \mathcal{O} \otimes \ZZ_2 = \mathcal{O}_2. $$
Let $z \in \WW$ be an element of norm $-1$.
Then the element $S = zT$ in $\mathcal{O}_2$ has the property $S^2 = 2$ and
$S a = \br{a} S$ for $a \in \WW$.  This makes explicit the presentation of
$\mathcal{O}_2$ in terms of $S$ and $\WW$
given in equation~(\ref{eq:Op}). 
\vspace{12pt}

\noindent
\emph{Claim 1:}
Let $a, b, \in \WW$ be such that $x = (1 + a) + bS \in \mathcal{O}_2$
has minimal polynomial $x^2 + 1$.
Then we must have $a \equiv 0 \pmod 2$ and $\nu_2(b) = 0$.
\vspace{12pt}

\noindent
\emph{Proof of Claim 1:}
In order for $x$ to have this minimal polynomial, we must have the following:
\begin{eqnarray}
\label{eq:itrace}
\Tr(x) &=& 2 + \Tr(a) = 0, \\
\label{eq:inorm}
N(x) &=& 1 + \Tr(a) + N(a) - 2N(b) = 1.
\end{eqnarray}
Reducing equation~(\ref{eq:itrace}) mod $4$, we see that 
$$ \Tr(a) \equiv 2 \pmod 4. $$
Substituting this into equation~(\ref{eq:inorm}) and reducing mod $2$, we see
that $N(a) \equiv 0 \pmod 2$.  Write $a = 2a'$.  Then we have $N(a) =
4N(a')$.  Therefore, when we reduce equation~(\ref{eq:inorm}) modulo 4,
we get $N(b) \equiv 1 \pmod 2$, and we conclude that $\nu_2(b) = 0$.
\vspace{12pt}

\noindent
\emph{Claim 2:}
Suppose $a, b \in \WW$ and $\alpha \in \ZZ_2$ are such that $x = (1 + a) + bS
\in \mathcal{O}_2$ has minimal polynomial $x^2 + \alpha x + 1$.
Then if $\alpha$ satisfies $\alpha \equiv 6 \pmod
{16}$, we must have $a \equiv 0 \pmod 2$ and $\nu_2(b) = 1$.
\vspace{12pt}

\noindent
\emph{Proof of Claim 2:}
Because $\alpha \equiv 6 \pmod{16}$, we can write $2 + \alpha = 4\gamma$,
where $\gamma \equiv 2 \pmod 4$.  

In order for $x$ to have this
minimal polynomial, we must have $\Tr(x) = -\alpha$ and $N(x) = 1$.
As a result, we have the following:
\begin{eqnarray}
\label{eq:1}
\Tr(a) &=& -4\gamma,\\
\label{eq:2}
N(a) &=& 2N(b) + 4\gamma.
\end{eqnarray}

Reducing equation~(\ref{eq:2}), we find that
\begin{eqnarray}
\label{eq:3}
N(a) &\equiv& 2N(b) \pmod 4.
\end{eqnarray}
This shows that $N(a) \equiv 0 \pmod 2$.  Writing $a = 2a'$ and
substituting back into equation~(\ref{eq:3}), we find that $N(b)
\equiv 0 \pmod 2$, so $b = 2b'$.

Re-expanding equations~(\ref{eq:1})
and
(\ref{eq:2}) gives the following:
\begin{eqnarray}
\label{eq:4}
\Tr(a') &=& -2\gamma,\\
\label{eq:5}
N(a') &=& 2N(b') + \gamma.
\end{eqnarray}

From equation~(\ref{eq:4}), we find $\Tr(a') \equiv 0 \pmod 2$.
Write $a' = a_1 + a_2\omega \in \WW$.
Because $\Tr(a') = 2a_1 - a_2 \equiv 0 \pmod 2$, we find $a_2 \equiv
0 \pmod 2$.  As a result we can write $a' = u + vd$ for 
$u,v \in \ZZ_2$ and $d = \sqrt{-3} =
2\omega+1$.  We then have $N(a') = u^2 + 3v^2$.  Therefore,
equation~(\ref{eq:5}) can be
reduced as follows:
\begin{eqnarray}
\label{eq:7}
u^2 - v^2 \equiv 2N(b') + 2 \pmod 4.
\end{eqnarray}

However, the equation $u^2 - v^2 \equiv 2 \pmod 4$ has no integer
solutions.  In order for equation~(\ref{eq:7}) to hold, we must have
$N(b') \equiv 1 \pmod 2$, or equivalently $\nu_2(b) = 1 + \nu_2(b') =
1$.
\vspace{12pt}

\noindent
\emph{Claim 3:}
Given $x \in \MSyL_2$, let $x' = \omega^2x\omega$.  Then $x'$ is in
$\MSyL_2$ and we have
$$
t_i(x') = \begin{cases}
t_i(x) & \text{$i$ even,} \\
\omega t_i(x) & \text{$i$ odd.}
\end{cases}
$$
This claim is immediate from the definition of the functions $t_i$ given in
Section~\ref{sec:density}.
\vspace{12pt}

We now complete the proof of Theorem~\ref{thm:Lambda} for the case $p=2$.
Consider the polynomial
$$ f(x) = x^2 + \frac{6}{\ell^4} x + 1 $$
with discriminant $\Delta = 4(9/\ell^8 - 1)$.  We have that $\Delta <
0$, so the polynomial $f$ is
irreducible over $\RR$, and $f$ is irreducible over $\ZZ_2$ because
$\nu_2(\Delta)$ is odd.  By Proposition~\ref{prop:Swan} and
Corollary~\ref{cor:minpoly} there
exists an element $y$ of $\mathcal{O}[1/\ell]$ so that $f(y) = 0$.  
Because $N(y) = 1$, the element $y$ lies in $\Gamma^1$.

In order to show that 
$\Lambda$ satisfies the hypotheses of Corollary~\ref{cor:density_2}, we claim
that the elements 
$$ i, k, y, y'=\omega^2 y \omega $$
lie in $\Lambda$, and their images under the homomorphism
$$ t_1 \oplus(t_3 + t_1t_2): \MSyL_2 \rightarrow \FF_4 \oplus \FF_4 $$
form an $\FF_2$-basis of $\FF_4\oplus \FF_4$.
Claims~1, 2, and 3 imply that these elements do lie in $\Lambda$, and the
functions $t_i$ evaluated on them satisfy
\begin{align*}
t_1(i) & \ne 0, \\
t_1(k) & = \omega t_1(i), \\
t_1(y) & = 0, \\
t_1(y') & = 0, \\
t_3(y) & \ne 0, \\
t_3(y') & = \omega t_3(y).
\end{align*}
These conditions are sufficient to conclude that their images give a basis.
Corollary~\ref{cor:density_2} now implies that $\Lambda$ is dense in
$\MSyL_2$.
\qed \end{pf*}

\begin{pf*}{Proof of Corollary~\ref{cor:Gamma1}.}
There is a short exact sequence
$$ 1 \rightarrow \MSyL_2 \rightarrow \MSL_2 \rightarrow C_{p+1}
\rightarrow 1,$$
where the cyclic group $C_{p+1}$ is the group of elements of 
$\FF_{p^2}^\times = (\mathcal{O}_p/(S))^\times$ of $\FF_p$-norm $1$.
It therefore suffices to show that we can lift the generator of $C_{p+1}$ to an
element of $\Gamma^1$.  Let $\br{y} \in \FF_{p^2}^\times$ be a generator of
the norm $1$ subgroup, with minimal polynomial
$$ \br{f}(x) = x^2 + ax + 1 $$
over $\FF_p$.
Let $\tilde a$ be an integer that reduces to $a$ modulo $p$, and 
define $\alpha = \tilde a/\ell^{m(p-1)}$
where $m$ is chosen
sufficiently large so that 
$\alpha^2 < 4$.
Then the polynomial
$$ f(x) = x^2 + \alpha x + 1 $$
is irreducible over $\QQ_p$ and $\RR$.
Just as in the proof of Theorem~\ref{thm:Lambda}, 
Proposition~\ref{prop:Swan} and Corollary~\ref{cor:minpoly} may be used to
show that there exists an element $y \in \mathcal{O}[1/\ell]$
so that $y$ reduces to a generator of $C_{p+1}$.
Because $y$ has norm $1$, it lies in $\Gamma^1$.
\qed \end{pf*}

\begin{pf*}{Proof of Theorem~\ref{thm:Gamma}.}
Assume that $p > 2$.
In light of Corollary~\ref{cor:Gamma1} and the short exact sequence
$$ 1 \rightarrow \MSL_2 \rightarrow \MS_2 \xrightarrow{N} \ZZ_p^\times 
\rightarrow 1,$$
we must show that there exists an element $x$ of $\Gamma$ so that $N(x)$ is
a topological generator of $\ZZ_p^\times$.
Because $\ell$ was assumed to be a topological generator, it suffices to show
there exists an $x$ so that $N(x) = \ell$.  
By
Proposition~\ref{prop:Kohel}, for $m$ sufficiently large there exists an
endomorphism $\alpha \in \End(C)$ of degree $\ell^{2m+1}$.  Then the
element $x = \ell^{-m}\alpha \in \Gamma$ has norm $\ell$.  
The argument for
$p=2$ is identical, except that we use the short exact sequence
$$ 1 \rightarrow \MSL_2 \rightarrow \td{\MS}_2 \xrightarrow{N}
\ZZ_2^\times/\{ \pm 1\} \rightarrow 1. \hfill \qed$$
\end{pf*}


\begin{thebibliography}{99}

\bibitem{Baker}
A. Baker,
Isogenies of supersingular elliptic curves over finite fields and
operations in elliptic cohomology.  Glasgow University Mathematics
Department preprint 98/39.

\bibitem{BehrensModular}
M. Behrens, 
A modular description of the $K(2)$-local sphere at the prime $3$.
To appear in Topology.

\bibitem{BehrensTree}
M. Behrens,
Isogenies of elliptic curves, Buildings, and the $K(2)$-local sphere.  
In preparation.

\bibitem{Deuring}
M. Deuring,
Die Typen der Multiplikatorenringe elliptischer Funktionenk\"orper. 
Abh. Math. Sem. Hansischen Univ. 14, (1941). 197--272.

\bibitem{DevinatzHopkins}
E.S. Devinatz, M.J. Hopkins,
Homotopy fixed point spectra for closed subgroups of the Morava 
stabilizer groups. Topology 43 (2004), no. 1, 1--47.

\bibitem{Giambalvo}
V. Giambalvo, 
Some tables for formal groups and $BP$. 
Geometric applications of homotopy theory 
(Proc. Conf., Evanston, Ill., 1977), II, pp. 169--176, 
Lecture Notes in Math., 658, Springer, Berlin, 1978.

\bibitem{GHMR}
P. Goerss, H.-W. Henn, M. Mahowald, C. Rezk,
A resolution of the $K(2)$-local sphere.  To appear in Ann. of Math.

\bibitem{GoerssHopkins}
P.G. Goerss, M.J. Hopkins, 
Moduli spaces of commutative ring spectra. 
Structured ring spectra, 151--200, London Math. Soc. Lecture Note Ser., 
315, Cambridge Univ. Press, Cambridge, 2004.

\bibitem{GMS}
V. Gorbounov, M. Mahowald, P. Symonds, 
Infinite subgroups of the Morava stabilizer groups. 
Topology 37 (1998), no. 6, 1371--1379.

\bibitem{Koch}
H. Koch,
Galoissche Theorie der $p$-Erweiterungen. 
Mit einem Geleitwort von I. R. \v Safarevi\v c. 
Springer-Verlag, Berlin-New York; 
VEB Deutscher Verlag der Wissenschaften, Berlin, 1970.

\bibitem{Kohelthesis}
D.R. Kohel,
Endomorphism rings of elliptic curves over finite fields.
Ph.D. Thesis, University of California at Berkeley.

\bibitem{Kohel}
D.R. Kohel,
Hecke module structure of quaternions. 
Class field theory---its centenary and prospect (Tokyo, 1998), 177--195, 
Adv. Stud. Pure Math., 30, Math. Soc. Japan, Tokyo, 2001.

\bibitem{Milne}
J.S.~Milne, Points on Shimura varieties mod $p$. 
Automorphic forms, representations and $L$-functions 
(Proc. Sympos. Pure Math., Oregon State Univ., Corvallis, Ore., 1977), 
Part 2, pp. 165--184, Proc. Sympos. Pure Math., XXXIII, 
Amer. Math. Soc., Providence, R.I., 1979.

\bibitem{Morava}
J. Morava,
Noetherian localisations of categories of cobordism comodules. 
Ann. of Math. (2) 121 (1985), no. 1, 1--39.

\bibitem{Pizer}
A. Pizer,
An algorithm for computing modular forms on $\Gamma \sb{0}(N)$. 
J. Algebra 64 (1980), no. 2, 340--390.

\bibitem{RavenelStructure}
D.C. Ravenel,
The structure of Morava stabilizer algebras. 
Invent. Math. 37 (1976), no. 2, 109--120. 

\bibitem{RavenelCohomology}
D.C. Ravenel,
The cohomology of the Morava stabilizer algebras. 
Math. Z. 152 (1977), no. 3, 287--297.

\bibitem{Ravenel}
D.C. Ravenel,
Complex cobordism and stable homotopy groups of spheres. Second edition.
AMS Chelsea Publishing, Amer. Math. Soc., Providence, RI, 2004.

\bibitem{Ravenelorange}
D.C. Ravenel,
Nilpotence and periodicity in stable homotopy theory. 
Appendix C by Jeff Smith. 
Annals of Mathematics Studies, 128. Princeton University Press, 
Princeton, NJ, 1992. 

\bibitem{RezkHMT}
C. Rezk,
Notes on the Hopkins-Miller theorem. 
Homotopy theory via algebraic geometry and group representations 
(Evanston, IL, 1997), 313--366, 
Contemp. Math., 220, Amer. Math. Soc., Providence, RI, 1998.

\bibitem{Riehm}
C. Riehm,
The norm $1$ group of a $\mathfrak{p}$-adic 
division algebra. Amer. J. Math. 92 1970 499--523. 

\bibitem{Serre}
J.P. Serre,
Cohomologie galoisienne. 
Cours au Coll\`ege de France, 1962-1963. Seconde \'edition. 
With a contribution by Jean-Louis Verdier. 
Lecture Notes in Mathematics, Vol. 5. Springer-Verlag, 
Berlin-Heidelberg-New York 1964.

\bibitem{Silverman}
J.H. Silverman,
The arithmetic of elliptic curves. 
Graduate Texts in Mathematics, 106. Springer-Verlag, New York, 1986.

\bibitem{Swan}
R.G. Swan,
$K$-theory of finite groups and orders. 
Lecture Notes in Mathematics, Vol. 149. Springer-Verlag, 
Berlin-New York, 1970.

\bibitem{Waterhouse-Milne}
W.C. Waterhouse, J.S. Milne,
Abelian varieties over finite fields.
Proc. Sympos. Pure Math. XX (1971), 53--64.

\end{thebibliography}
\end{document}